\documentclass[11pt]{article}
\usepackage[utf8]{inputenc}
\usepackage[english]{babel}
\usepackage{amsmath,amsfonts,amscd,amssymb}
\usepackage[all]{xy}
\usepackage{makeidx}
\usepackage[paper=a4paper,left=30mm,right=20mm,top=25mm,bottom=30mm]{geometry}

\newtheorem{thm}{Theorem}[section]
\newtheorem{cor}[thm]{Corollary}
\newtheorem{lem}[thm]{Lemma}
\newtheorem{prop}[thm]{Proposition}
\newtheorem{deft}[thm]{Definition}
\newtheorem{rek}[thm]{Remark}
\newtheorem{conj}[thm]{Conjecture}
\let\noi=\noindent

\def\N{\mathbb{N}} 

\makeindex
\title{A description of and an upper bound on the set of bad primes in the study of the Casas-Alvero Conjecture}
\author{Daniel Schaub, Univ Angers, CNRS, LAREMA, SFR MATHSTIC\\
F-49000 Angers, France
\\email: daniel.schaub@univ-angers.fr \and
Mark Spivakovsky, Univ Paul Sabatier, CNRS, IMT UMR 5219\\
F-31062 Toulouse, France and\\
CNRS, LaSol UMI 2001, UNAM.}

\begin{document}
\maketitle

\begin{abstract} The Casas--Alvero conjecture predicts that every univariate polynomial over a field of characteristic zero having a common factor with each of its derivatives $H_i(f)$ is a power of a linear polynomial. One approach to proving the conjecture is to first prove it for polynomials of some small degree $n$, compile a list of bad primes for that degree (namely, those primes $p$ for which the conjecture fails in degree $n$ and characteristic $p$) and then deduce the conjecture for all degrees of the form $np^\ell$, $\ell\in\N$, where $p$ is a good prime for $n$. In this paper we give an explicit description of the set of bad primes in any given degree $n$. In particular, we show that if the conjecture holds in degree $n$ then the bad primes for $n$ are bounded above by $\binom{\frac{n^2-n}2}{n-2}!\prod\limits_{i=1}^{n-1}
\binom{i+n-2}{n-2}^{\binom{d-i+n-2}{n-2}}$.
\end{abstract}

\section{Introduction}

In the year 2001 Eduardo Casas--Alvero published a paper on higher order polar germs of plane curve singularities \cite{C}. His work on polar germs inspired him to make the following conjecture.

Let $K$ be a field, $f\in K[x]$ a non-constant monic univariate polynomial, $n:=\deg(f)$:
\[
f=x^n+a_1x^{n-1}+\dots+a_n.
\]
Let
$$
H_i(f)= \binom{n}{i}x^{n-i} + \binom{n-1}{i}a_1x^{n-i-1} + \cdots + \binom{i}{i}a_{n-i}
$$
be the $i$-th Hasse derivative of $f$.

\begin{deft}
The polynomial $f$ is said to be a {\bf Casas--Alvero polynomial} if for each $i\in\{1,\ldots,n-1\}$ it has a non-constant common factor with its $i$-th Hasse derivative $H_i(f)$.
\end{deft}

\begin{conj} {\bf (Casas--Alvero)}
Assume that $\text{char}\ K=0$. If $f \in K[x]$ is a Casas-Alvero polynomial of degree $n$, then there exists $b\in K$ such that $f(x) =(x-b)^n$.
\end{conj}

If $\text{char}\ K=p >0$, the conjecture is false in general. The simplest counterexample is the polynomial $f(x)= x^{p+1}-x^p$.

\begin{rek}
The following fact is known and easy to prove. If the Casas--Alvero conjecture holds for the algebraic closure $\bar K$ of $K$ then it also holds for $K$ (the converse is not established, to our knowledge). Therefore, from now on we will assume that $K$ is algebraically closed.
\end{rek}

We will write CA$_{n,p}$ for the statement ``The Casas--Alvero conjecture holds for polynomials of degree $n$ over algebraically closed fields of characteristic $p$''.

The following equivalences are known for each $n$ (\cite{DdJ}, \cite{graf}):

CA$_{n,0}$ holds $\iff$ CA$_{n,p}$ holds for some prime n

\begin{deft}
A prime number $p$ is said to be a \textbf{bad prime for} $n$ if CA$_{n,p}$ is false. If $p$ is not a bad prime for $n$, we will say that $p$ is a \textbf{good prime for} $n$.
\end{deft}

\begin{prop} (\cite{graf}, Propositions 2.2 and 2.6)
Take a strictly positive integer $n$, a prime number $p$ and a
non-negative integer $\ell$. Assume that CA$_{n,p}$ holds. Then so do CA$_{np^\ell,p}$ and CA$_{np^\ell,0}$.
 \end{prop}

\noi This result suggests the following general approach to the problem:

(1) prove the conjecture for a small number $n$;

(2) compile lists of good and bad primes for $n$;

(3) conclude that CA$_{np^\ell,0}$ holds for all the primes $p$ that are known to be good for $n$.
\medskip

In particular, this shows the importance of knowing which primes are good or bad for a given degree $n$.

The above approach has been carried out up to $n \le 7$ (\cite{CLO}, \cite{C-S}, \cite{DdJ}, \cite{frutos}, \cite{graf}). Some integers cannot be written in the form $np^\ell$ where $p$ is a good prime for $n$, for example,
$$
12=2^2 \cdot3,\ 20=2^2\cdot5,\ 24=2^3\cdot3,\ 28=2^2\cdot7,\ 30=2\cdot3\cdot5,\ 36=2^2\cdot3^2,\
40=2^3\cdot5,\ldots
$$
CA$_{12,0}$ has been proved in \cite{CLO} with the aid of a computer, by using a very clever strategy to cut down the computation of resultants and Gröbner bases. Thus the smallest degree $n$ for which CA$_{n,0}$ is not known is $n=20$.
\medskip

The purpose of this paper is to give an explicit description of the set of bad primes in any given degree $n$. In particular, we obtain an explicit upper bound on bad primes for $n$, assuming that $C_{n,0}$ holds. These results are based on recent work of Soham Ghosh \cite{Go}.

\medskip

\noindent{\bf Notation:} For $j\in\{1,\dots,n-1\}$, let the involution
\begin{equation}
\Phi_j: K[x_1,\ldots,x_{n-1}] \to K[x_1,\ldots,x_{n-1}],
\end{equation}
be defined by
\begin{equation}
\Phi_j(x_i)  =x_i-x_j \text{ for } i \neq j \text{ and } \Phi_j(x_j)=-x_j.
\end{equation}
Let
\begin{equation}
\Phi_n: K[x_1,\ldots,x_{n-1}] \to K[x_1,\ldots,x_{n-1}],
\end{equation}
denote the identity map.

Let $\sigma_i(x_1,\ldots,x_{n-1})$ denote the $i$-th elementary symmetric function of $x_1,\ldots,x_{n-1}$.

Let $\mathcal T=\{1,\dots,n\}^{n-1}$; the set $\mathcal T$ is the collection of all the $(n-1)$-tuples of the form
$(j_1,\dots,j_{n-1})$, where $j_1,\ldots,j_{n-1}\in\{1,\dots,n\}$.
\medskip

\noindent{\bf Notation.} Given a fixed choice of
$T=(j_1,\dots j_{n-1})$, for $i\in\{1,\dots,n-1\}$ we will denote by $G_{T,i}$ the homogeneous polynomial
$\Phi_{j_i}(\sigma_i(x_1,\ldots,x_{n-1}))$.
\medskip

In his fundamental preprint \cite{Go}, Soham Ghosh showed that the Casas--Alvero conjecture in degree $n$ (over any field, regardless of characteristic) is equivalent to the following statement.

\begin{conj}\label{mainthm} (\cite{Go}, Proposition 5.2) For every choice of $T=(j_1,\ldots,j_{n-1})\in\mathcal T$, the sequence of homogeneous polynomials
\begin{equation}\label{eq:suite}
    \left(G_{T,1},\ldots,G_{T,{n-1}}\right)
\end{equation}
forms a regular sequence in $K[x_1,\ldots,x_{n-1}]$.
\end{conj}

Since the polynomial ring $K[x_1,\ldots,x_{n-1}]$ is
Cohen--Macaulay, Conjecture \ref{mainthm} is equivalent to saying that $\text{ht}(G_{T,1},\ldots,G_{T,{n-1}})=n-1$ and thus also to

\begin{conj}\label{radical}
We have
\begin{equation}
\sqrt{G_{T,1}, \ldots, G_{T,{n-1}}}=
(x_1,\dots x_{n-1}).\label{eq:regseq}
\end{equation}
\end{conj}

\section{Macaulay's Theorem}

We recall (a part of) Macaulay's celebrated theorem from 1916.

Let $x_1,\ldots,x_n$ be independent variables, $f_1,\dots,f_n\in K[x_1,\ldots,x_n]$ homogeneous polynomials and let $d_i=\deg\ f_i$ denote the total degree of $f_i$. Let
$\mathfrak m:=(x_1,\ldots,x_n)$. Finally, put $d=\sum\limits_{k=1}^nd_k-n+1$.

\begin{thm}\cite{Mac}\label{macaulay} The following statements are equivalent:
\begin{enumerate}
    \item[(1)] $\sqrt{(f_1,\dots,f_n)}=\mathfrak m$
    \item[(2)] $\mathfrak m^d\subset(f_1,\dots,f_n)$.
\end{enumerate}

\end{thm}

\section{A description of and an upper bound on the set of bad primes}

In this section we state and prove our main results.

Let $x=(x_1,\dots,x_{n-1})$. We will use multi-index notation: $x^\alpha$ will stand for $\prod\limits_{k=1}^{n-1}x_k^{\alpha_k}$ and $|\alpha|$ for $\sum\limits_{k=1}^{n-1}\alpha_k$.

We apply Macaulay's Theorem to the polynomials
$G_{T,},\dots,G_{T,{n-1}}\in K[x]$.

We have $\deg\ G_{T,i}=i$ for $i\in\{1,\dots,n-1\}$.

Let $d=\sum\limits_{i=1}^{n-1}\deg\ G_{T,i}-(n-2)=
1+2+\dots+(n-1)-(n-2)=\frac{n^2-3n+4}2$.

Let $C$ denote the binomial coefficient
$\binom{\frac{n^2-n}2}{n-2}$; it is the number of monomials of degree $d=\frac{n^2-3n+4}2$ in $n-1$ variables.

Let $S_{T,i}=\left\{\left.G_{T,i}x^\alpha\ \right|\
|\alpha|=d-i\right\}$, $i\in\{1,\dots,n-1\}$ and
\[
S_T:=\bigcup\limits_{i=1}^{n-1}S_{T,i}.
\]
We have $|S_{T,i}|=\binom{d-i+n-2}{n-2}$.

Let $D:=|S_T|=
\sum\limits_{i=1}^{n-1}|S_{T,i}|=
\sum\limits_{i=1}^{n-1}\binom{d-i+n-2}{n-2}$; we have $D\ge C$ (in fact, this inequality is strict whenever $n>2$).

Consider the $C$-dimensional $K$-vector space $V$, generated by all the monomials in $x$ of degree $d$; we have
$S_T\subset V$.

Let $M_T$ denote the $D\times C$ matrix formed by the row vectors $(v)_{v\in S_T}$. Let $J_T$ be the greatest common divisor of all the $C\times C$ minors of $M_T$.

\begin{thm}\label{badprimes}
A prime number $p$ is a bad prime for $n$ if and only if $p\ |\ J_T$ for some $T\in\mathcal T$ (equivalently, if and only if $p\ |\
{\rm lcm}(J_T)_{T\in\mathcal T}$).
\end{thm}
{\it Proof.} Fix a $T\in\mathcal T$. By Theorem \ref{macaulay},  (\ref{eq:regseq}) is equivalent to $(x)^d\subset(G_{T1},\dots,G_{T,{n-1}})$. And this is true if and only if $V\subset(G_{T1},\dots,G_{T,{n-1}})$. This inclusion is true if and only if the rank of the matrix $M_T$ is maximal, ie. $\text{rk}\ M_T=C$. or in other words, if and only $M_T$ has a non-degenerate $C\times C$ minor.

Therefore Conjecture \ref{mainthm} fails in degree $n$ and characteristic $p$ if and only if $p\ |\ J_T$ for some $T\in\mathcal T$. By \cite{Go}, Proposition 5.2, the failure of Conjecture \ref{mainthm} in degree $n$ and characteristic $p$ is equivalent to $p$ being a bad prime for $n$.
$\Box$
\medskip

\begin{cor}{\ }
If C$_{n,0}$ holds but $p$ is a bad prime for $n$ then
\begin{equation}
p<C!\prod\limits_{i=1}^{n-1}
\binom{i+n-2}{n-2}^{\binom{d-i+n-2}{n-2}}.\label{eq:upperbound}
\end{equation}
\end{cor}
\noi{\it Proof.} The corollary follows from Theorem \ref{badprimes} and the following lemma.

\begin{lem}
Fix a $T\in\mathcal T$ and let $A$ be a $C\times C$ minor of $M_T$. Then
\[
|A|\le C!\prod\limits_{i=1}^{n-1}
\binom{i+n-2}{n-2}^{\binom{d-i+n-2}{n-2}}.
\]
\end{lem}

\noindent{\it Proof of the lemma.} Write $G_{T,i}$ as a sum of (possibly repeated) monomials, each with coefficient 1. The monomial $x_{j_i}^i$ is repeated $\binom{i+n-2}{n-2}$ times, more than any other monomial. Therefore, once we group the like terms together,  the greatest absolute value of a coefficient of a monomial in $G_{Ti}$ is $\binom{i+n-2}{n-2}$.

When we write the minor $A$ as a sum of $C!$ terms, each term divides an integer of the form $\prod\limits_{i=1}^{n-1}
\prod\limits_{j=1}^{\binom{i+n-2}{n-2}}a_{ij}$, where for all the pairs $(i,j)$ we have $|a_{ij}|\le\binom{i+n-2}{n-2}$. This proves the lemma and, with it, the corollary.
$\Box$

\begin{rek}
The upper bound \eqref{eq:upperbound} can be vastly improved as follows. Let the notation be as above. The product $\prod\limits_{i=1}^{n-1}
\binom{i+n-2}{n-2}^{\binom{d-i+n-2}{n-2}}$ has a total of $D$
(not necessarily distinct) terms of the form $\binom{i+n-2}{n-2}$. We have $\binom{i+n-2}{n-2}<\binom{i'+n-2}{n-2}$ whenever $i<i'$. Write $\prod\limits_{i=1}^{n-1}
\binom{i+n-2}{n-2}^{\binom{d-i+n-2}{n-2}}=\prod\limits_{k=1}^Db_k$ with the sequence $b_k$ non-strictly increasing. Then
\begin{equation}
p<C!\prod\limits_{D-C+1}^Db_k\label{eq:improvement}
\end{equation}
The proof is the same as in the corollary. The reason we did not state the corollary in this form in the first place is that we could not find an explicit, closed form for the integers $b_k$ appearing in \eqref{eq:improvement}. For example, what is the value of $i$ such that $b_{D-C+1}=\binom{i+n-2}{n-2}$ ?
\end{rek}

\end{document}